\documentstyle{amsart}
\def\dj{d\kern-.30em\raise1.25ex\vbox{\hrule width .3em height .03em}}
\def\Dj{D\kern-.75em\raise0.75ex\vbox{\hrule width .3em height .03em}
\kern.03em}
\newsymbol\rest 1316
\def\sstar{{\raise0.2ex\hbox{$\scriptstyle\star$}}}
\newcommand{\hor}{\mbox{\family{euf}\shape{n}\selectfont hor}}
\newcommand{\SO}{\mbox{\shape{u}\selectfont SO}}
\newcommand{\so}{\mbox{\family{euf}\shape{n}\selectfont so}}
\newcommand{\s}{\sigma}
\newcommand{\e}{e}
\newcommand{\r}{\rho}
\newcommand{\cl}{C\!\ell}
\newcommand{\BC}{\Bbb{C}}
\newcommand{\Spin}{\mbox{S}}
\newcommand{\JJ}{\sqcup}
\newcommand{\Sum}{\displaystyle{\sum}}
\newcommand{\Int}{\int}
\def\intM{\int_M^\uparrow}
\newcommand{\ad}{\mbox{\shape{n}\selectfont ad}}
\newcommand{\id}{\mbox{\shape{n}\selectfont id}}
\newcommand{\hM}{{\star}_M}
\newcommand{\hE}{{\star}_E}
\newcommand{\nE}{\nabla_{\!E}}
\newcommand{\dM}{d_M}
\newcommand{\lM}{\Delta_M}
\newcommand{\lE}{\Delta_E}
\newcommand{\WM}{\Omega_M}
\newcommand{\EW}{{\cal E}_\Omega}
\renewcommand{\S}{\Sigma}
\newtheorem{lem}{Lemma}
\numberwithin{equation}{section}
\begin{document}
\title[Spinor Structures]{Classical Spinor
Structures on Quantum Spaces}
\author{Mi\'cO \Dj ur\Dj evi\'c}
\address{Instituto   de   Matematicas,
UNAM, Area de la Investigacion Cientifica,
Circuito Exterior, Ciudad Universitaria,
M\'exico DF, CP 04510, M\'EXICO\newline
\indent {\it Written In}\newline
\indent
Centro de Investigaciones Teoricas, UNAM,
Facultad de Estudios Superiores Cuautitlan, Cuautitlan Izcalli,
M\'EXICO}
\maketitle
\begin{abstract}
A noncommutative-geometric generalization of the classical concept
of spinor structure is presented. This is done in the framework of
the formalism of quantum principal bundles. In particular, analogs
of  the  Dirac operator and the Laplacian are introduced and
analyzed.  A  general construction of examples of quantum spaces
with a spinor structure is presented.
\end{abstract}
\section{Introduction}
\renewcommand{\thepage}{}
Classical geometry of spinor structures is based on two
extremely reach and effective mathematical concepts. As first,  there
is a ``local'',  purely  algebraic,  concept  of  spinors,  which
includes Clifford algebras and spin groups. The second one is a
``global'' concept of principal bundle, specified as a covering
spin-bundle of the bundle of oriented  orthonormal frames (of the base
space).

In this paper basic  concepts  of  classical spinor geometry
will  be  translated  into  an  appropriate  quantum  context,  in
accordance with general principles of
non-commutative differential geometry \cite{C1,C2}.

The starting idea is that in formulating the concept of a spinor
structure on a quantum  space,  only  ``global''  aspects  of  the
classical theory should be quantized. In formalizing this we shall
use  a  general  non-commutative-geometric  theory  of   principal
bundles, presented in \cite{D1,D2,D3}.
Spinor structures on quantum
spaces will be represented by certain quantum  principal  bundles,
however the structure group will be a  classical  spin  group.

In this sense, ``local aspects'' of the classical theory will be
left unchanged. The reason for this is  that  introducing  in  the
theory a quantum deformation of the spin group would cause the
lack of internal symmetry,  because
of the inherent geometrical inhomogeneity of (truly) quantum
groups \cite{D1}. Such a philosophy opens the possibility to generalize,
in a more or less straightforward way, fundamental concepts and
constructions of the classical theory (presented in \cite{Cru},
for example).

The paper is organized as follows. In the next section the concept
of a quantum spin manifold will be introduced, in the
framework of the formalism of framed  quantum  principal   bundles
\cite{frm}.
Informally speaking, a spinor structure on  a  quantum  space  $M$
will be represented by a principal spin-bundle $P$ over $M$,
endowed with an additional structure,  which  expresses  the  idea
that $M$ is an ``oriented Riemannian manifold'' such that $P$
``covers'' the bundle of oriented  orthonormal
frames  for  $M$ (however,  this   bundle  will  figure  only
implicitly in the theory).

Further, we shall introduce the analog of the covariant derivative
of the Levi-Civita  connection  (operating on  $P$),  and  sketch  a
construction of the canonical differential calculus on  $P$.  Then
counterparts of the Hodge *-operator, Laplace  operator  and
invariant integration on $P$ will  be  constructed,  in  a  direct
analogy with classical differential geometry.

Section 3 deals with ``spinor fields''.  Starting  from
$P$  we  shall  introduce  a  space  $\cal  E$,  interpretable  as
consisting of ``smooth sections''  of  the corresponding
``associated'' spinor bundle $E$.
 A  quantum  counterpart  of  the  Dirac
operator will be then introduced and analyzed. In particular,
it turns out that the difference between the square of  the  Dirac
operator and the Laplacian (acting in $E$) is proportional to
the (scalar) curvature of  the  Levi-Civita  connection,  as
in  classical geometry \cite{Li}. We  shall  also define the
space of ``$E$-valued'' differential forms on $M$, and the
quantum counterpart of the Clifford bundle.

In  Section  4  the  theory  is  illustrated  in  a  simple,   but
sufficiently reach class  of  examples  of  inherently  ``curved''
quantum spin manifolds.

Finally, in Section 5 some concluding remarks are made.
\section{General Considerations}
\renewcommand{\thepage}{\arabic{page}}
In this section we shall introduce and analyze  spinor  structures
on quantum spaces. As first, let us fix the notation.
The Clifford algebra associated
with the complex Euclidean space $\Bbb{C}_n$ will be denoted by $\cl_n$.

 Since in  all  formal  considerations
quantum spaces are represented  ``dually''  by  the  corresponding
functional algebras, we shall also use this description for the
classical  groups  $\Spin(n)$ (the spin groups) and  $\SO(n)$,
in parallel  with  the
standard one. A (commutative)
*-algebra of polynomial functions on $\Spin (n)$ will be denoted by
${\cal S}_n$. The group structure on $\Spin (n)$ induces a Hopf  algebra
structure on ${\cal S}_n$, which is determined by the coproduct
$\phi\colon
{\cal S}_n\rightarrow{\cal S}_n\otimes{\cal S}_n$, counit
$\epsilon\colon{
\cal S}_n
\rightarrow\BC$ and the antipode
$\kappa\colon{\cal   S}_n\rightarrow{\cal   S}_n$   (we   follow    the
notation of \cite{W}, although in the classical context). In the
``dual'' picture,
the  elements  of  $\Spin  (n)$  are  naturally  interpretable  as
*-characters    $g\colon{\cal    S}_n\rightarrow\BC$ (non-trivial
hermitian multiplicative linear functionals).
Furthermore $\SO(n)$  is  represented
by a *-Hopf  subalgebra  ${\cal A}_n$  of  ${\cal S}_n$,
consisting  of
polynomial functions on $\Spin(n)$ invariant under the map
$g\mapsto -g$. At the
level of spaces, the inclusion ${\cal A}_n\hookrightarrow
{\cal S}_n$ becomes
the canonical (universal covering) epimorphism $\amalg\colon \Spin
(n)\rightarrow  \SO(n)$. We shall denote by $u$
the canonical
representation of $\Spin (n)$  in
$\BC_n$ (obtained by composing $\amalg$ and the standard
representation of $\SO(n)$ in $\BC_n$). Explicitly,
$$ u_g(x)=gxg^{-1}$$
for each  $x\in\BC_n$  and  $g\in \Spin (n)$  (and  $\Spin (n)$  and
$\BC_n$  are  embedded   in   $\cl_n$).   Equivalently,   $u$   is
understandable as a map
(a right ${\cal S}_n$-comodule structure on $\BC_n$) of the form
$u\colon\BC_n\rightarrow\BC_n\otimes   {\cal    S}_n$
so that
$u_g=(\id\otimes g)u$. Explicitly,
$$u(\e_i)=\sum_{j=1}^n \e_j\otimes u_{ji}$$
where $u_{ji}\in{\cal S}_n$ are matrix elements of $u$
(actually $u_{ij}$ generate the *-algebra ${\cal A}_n$).

Let  $M$  be  a  quantum  space represented by a
non-commutative *-algebra
${\cal V}$. The elements of $\cal V$ are
interpretable as ``smooth functions'' on $M$.  Let  $P=({\cal
B},i,F)$ be a (quantum) principal $\Spin (n)$-bundle
over $M$. Here,
${\cal B}$ is a *-algebra consisting of  appropriate
``functions'' on the quantum space $P$, while $i\colon{\cal V}
\rightarrow{\cal B}$ is a *-monomorphism playing the role of ``the
dualized projection'' of $P$ on $M$. Finally,
$F\colon{\cal B}\rightarrow{\cal B}\otimes{\cal S}_n$  is  a
*-homomorphism playing the role of the dualized right  action of
$\Spin (n)$ on $P$. This interpretation is  formalized
in the following conditions
\begin{align}
(\id\otimes\phi)F&=(F\otimes \id)F\\
(\id\otimes \epsilon)F&=\id.
\end{align}
Speaking geometrically, the space $M$ can be identified,  via  the
``projection map'' with the ``orbit space'' corresponding
to the action $F$.  Formally,  it  means  that  $i({\cal  V})$  is
consisting precisely of those elements $b\in{\cal B}$ which are
$F$-invariant, in  the  sense  that  $F(b)=b\otimes  1$.  In  what
follows, we shall identify the elements of  $\cal  V$  with  their
images in $i(\cal V)$.

The action $F$ is ``free'' in the sense
that for  each  $a\in{\cal
S}_n$ there exist elements $q_i,b_i\in{\cal B}$ satisfying
\begin{equation}\label{free}
\sum_iq_iF(b_i)=1\otimes a.
\end{equation}
The actual action of elements of $\Spin (n)$ on ${\cal B}$  is
described by *-automorphisms $$F_g=(\id\otimes g)F.$$

By definition \cite{frm} a {\it frame structure} on the bundle
$P$  (relative
to $u$) is every
$n$-tuple   $\tau=(\partial_1,\ldots,\partial_n)$   of   hermitian
derivations    $\partial_i\colon{\cal    V}\rightarrow{\cal    B}$
satisfying
$$ F\partial_i(f)=\sum_{j=1}^n\partial_j(f)\otimes u_{ji} $$
for each $f\in{\cal V}$ and $i\in\{1,\ldots,n\}$ and such that
there exist elements $b_{i\alpha}\in{\cal B}$ and
$v_{i\alpha}\in{\cal V}$ with the property
$$ \sum_{\alpha}b_{i\alpha}\partial_j(v_{i\alpha})=\delta_{ij}1 $$
for each $i,j\in\{1,\ldots,n\}$.

A frame structure $\tau$ is  called  {\it  integrable}  iff there
exists a  system  $\widehat{\tau}=(X_1,\ldots,X_n)$  of
hermitian  derivations  $X_i\colon{\cal   B}\rightarrow{\cal   B}$
which extend derivations $\partial_i$ and such  that  the  following
identities hold
\begin{align*}
0&=X_i\partial_j-X_j\partial_i\\
FX_i(b)&=\sum_{j=1}^n\sum_k X_j(b_k)\otimes u_{ji}c_k,
\end{align*}
where $\Sum_k b_k\otimes c_k=F(b)$.

Speaking geometrically, frame structures formalize  the
idea that $M$ is an ``oriented Riemannian manifold''. Furthermore,
endowed with a frame structure $\tau$, the bundle $P$
becomes a ``covering bundle'' of the bundle $P^\sstar$ of  oriented
orthonormal frames for $M$.

Starting from  an integrable frame structure
$\tau$  it  is  possible  to  construct
\cite{frm} the  whole
differential calculus on the bundle
(including the calculus on $M$). The construction has several
steps. As first, a graded *-algebra $\hor_P$ representing
horizontal forms
on $P$ can be defined as
$$ \hor_P={\cal B}\otimes \BC_n^\wedge$$
where $(\,)^\wedge$ denotes the corresponding external algebra.
The *-structure on
$\BC_n^\wedge$ is specified by  $\e_i^*=\e_i$  and  extended  by
antilinearity and multiplicativity on the whole $\BC_n^\wedge$.
Then the formula
\begin{equation}
\nabla(b\otimes\vartheta)=\sum_{k=1}^n X_k(b)\otimes
(\e_k\wedge\vartheta)
\end{equation}
consistently  defines  a  hermitian   first-order   antiderivation
$\nabla$ on $\hor_P$. This map plays the  role  of  the  covariant
derivative (induced by  the
Levi-Civita connection).

There exists a natural action $F^\sstar\colon\hor_P\rightarrow
\hor_P\otimes{\cal S}_n$ of $\Spin (n)$
on $\hor_P$. It is defined by
\begin{equation}
F^\sstar(b\otimes\vartheta)=F(b)u^\wedge(\vartheta),
\end{equation}
where
$u^\wedge\colon\BC_n^\wedge\rightarrow\BC_n^\wedge\otimes
{\cal S}_n$ is the representation induced by $u$
(and $\cal B$, $\Bbb{C}_n^\wedge$ are understood as
subalgebras of $\hor_P$).
The map $F^\sstar$ is a *-homomorphism and
\begin{gather}
(\id\otimes\phi)F^\sstar=(F^\sstar\otimes \id)F^\sstar\\
(\id\otimes \epsilon)F^\sstar=\id\\
(\nabla\otimes \id)F^\sstar=F^\sstar\nabla\label{nab-cov}.
\end{gather}

Moreover, there  exists  the  unique
map $\varrho\colon{\cal S}_n\rightarrow\hor_P$ such that
\begin{equation}\label{nab2}
\nabla^2(\varphi)=-\sum_k \varphi_k \varrho(c_k)
\end{equation}
where $F^\sstar(\varphi)=\Sum_k\varphi_k\otimes c_k$.
The map $\varrho$ is interpretable as the corresponding ``curvature
tensor''. We have
\begin{align}
\varrho(1)&=0\label{varrho1}\\
F^\sstar\varrho&=(\varrho\otimes \id)\ad
\end{align}
where $\ad\colon{\cal S}_n\rightarrow{\cal  S}_n\otimes{\cal  S}_n$
is the dualized adjoint action of $\Spin (n)$ on itself. In what follows
it will be assumed that $\varrho$ is {\it classical}, in the sense that
\begin{equation}
\varrho(ab)=\epsilon(a)\varrho(b)+
\epsilon(b)\varrho(a)\label{r-der}
\end{equation}
for each $a,b\in\cal{A}$. Then the following properties hold
\begin{gather}
\varphi\varrho(a)=\varrho(a)\varphi\\
\varrho(a^*)=\varrho(a)^*\\
\varrho(\kappa(a))=-\varrho(a).
\end{gather}

Let $\WM$ be the *-subalgebra of $\hor_P$ consisting of
$F^\sstar$-invariant elements. According to \eqref{nab-cov}
we have  $\nabla(\WM)\subseteq\WM$.  Let
$\dM\colon\WM\rightarrow\WM$  be  the
restriction of $\nabla$ on $\WM$. This  map  is  a  hermitian
first-order antiderivation. Furthermore,
$$d_M^2=0$$
(according to \eqref{nab2} and \eqref{varrho1}).
The elements of $\WM$ are interpretable as ``differential
forms''on $M$. The map $\dM$ plays the role of the exterior
derivative.

A  graded-differential   *-algebra   $\Omega_P$   describing   the
whole differential calculus on the bundle can be now constructed
as follows.
At the level of graded *-algebras,
$$\Omega_P=\hor_P\mathbin{\widehat{\otimes}}L^\wedge $$
where $\widehat{\otimes}$ denotes the graded tensor product.
In the above formula
$L=\so(n)^*$  (where  $\so(n)$  is
the Lie algebra of complex  antisymmetric  $n\times  n$  matrices,
naturally understood as the Lie algebra of
$\SO(n)$) and
$L^\wedge$ is endowed with
its natural graded-differential *-algebra structure.
The differential
$d\colon L^\wedge\rightarrow L^\wedge$  is  explicitly  given  by  the
Mauer-Cartan formula (on the first-order elements and extended  on
the whole $L^\wedge$ by the graded Leibniz rule).
Also,   $\hor_P$   and   $L^\wedge$   are
naturally understandable as subalgebras of $\Omega_P$.

The differential
$d^\wedge\colon\Omega_P
\rightarrow\Omega_P$ is specified by
\begin{align}
d^\wedge(\varphi)
&=\nabla(\varphi)+(-)^{\partial\varphi}
\sum_k\varphi_k\pi(c_k)\label{d-hor}\\
d^\wedge(\vartheta)&=R(\vartheta)+d(\vartheta)\label{d-L}
\end{align}
where $\vartheta\in L$ and $\varphi\in\hor_P$, while
$F^\sstar(\varphi)=\Sum_k\varphi_k\otimes c_k$. Next,
$\pi\colon{\cal S}_n\rightarrow L$ is a
natural projection, given by $\pi(a)(X)=(Xa)_\epsilon$ (here elements
$X\in \so(n)$ are understood as  left-invariant  vector  fields  on
$\Spin(n)$). Finally $R\colon{\cal
A}\rightarrow\hor_P$ is a map given by $R\pi=\varrho$
(consistency of this definition is a consequence
of \eqref{r-der}). Using the above
formulas and the graded Leibniz rule the map  $d^\wedge$  can  be
consistently extended to the whole $\Omega_P$.

Let $\varpi^\wedge\colon L^\wedge\rightarrow L^\wedge\otimes
{\cal S}_n$ be the
*-homomorphism   extending   the   dualized    coadjoint    action
$\varpi\colon L\rightarrow L\otimes {\cal S}_n$ (explicitly given by
$\varpi\pi=(\pi\otimes\id)\ad$).
The formula
\begin{equation}
F^\wedge(\varphi\otimes\vartheta)=F^\sstar(\varphi)
\varpi^\wedge(\vartheta)
\end{equation}
defines a *-homomorphism
$F^\wedge\colon\Omega_P\rightarrow\Omega_P
\otimes  {\cal  S}_n$  extending  the  action $F^\sstar$. We have
\begin{gather}
F^\wedge d^\wedge=(d^\wedge\otimes \id)F^\wedge\\
(F^\wedge\otimes \id)F^\wedge=(\id\otimes \phi)F^\wedge\\
\id=(\id\otimes \epsilon)F^\wedge
\end{gather}
The construction of $\Omega_P$ can be viewed as a special
case of a general construction of differential calculus
presented in \cite{D2}--Section 6.

A map $\omega\colon L\rightarrow\Omega_P$ defined by
$$\omega(\vartheta)=1\otimes\vartheta$$
is a regular and multiplicative connection on $P$.
It is a  counterpart
of  the  Levi-Civita  connection
(viewed as a connection on the spin-bundle)
in classical geometry. Moreover the
map $R$  corresponds to the {\it curvature } of
$\omega$, and $\nabla$ can be viewed as the {\it covariant
derivative} associated to $\omega$.

The quantum analog of the  Hodge  *-operator can  be  introduced
in a direct analogy with classical geometry. The formula
\begin{equation} {\star}\bigl(b\otimes (\e_{i_1}\wedge\ldots
\wedge\e_{i_k})\bigr)=b\otimes (\e_{j_1}\wedge\ldots \wedge\e_{j_l})
\end{equation}
where $1\leq i_1<\ldots<i_k\leq n$ and $1\leq j_1<\ldots<
j_l \leq n$, while $k+l=n$ and
$\e_{i_1}\wedge\ldots \wedge\e_{i_k}\wedge
\e_{j_1}\wedge\dots\wedge\e_{j_l}
=\e_1\wedge\ldots\wedge\e_n$, consistently determines a linear map
${\star}\colon\hor_P\rightarrow\hor_P$.

Let $\nabla^\dagger\colon\hor_P\rightarrow\hor_P$ be a linear
map defined by
\begin{equation}\label{div}
\nabla^\dagger(\varphi)=-(-)^{n-nk}
{\star}\nabla {\star}(\varphi)
\end{equation}
where $\varphi\in\hor_P^k$.

By construction maps ${\star}$ and $\nabla^\dagger$ are
right-covariant, in the sense that
\begin{align}
F^\sstar{\star}&=({\star}\otimes \id)F^\sstar\label{inv-*}\\
F^\sstar\nabla^\dagger&=(\nabla^\dagger\otimes\id)F^\sstar.
\end{align}
In particular ${\star}$ and $\nabla^\dagger$ are reduced in $\WM$.

Let us denote by ${\star}_M,d_M^\dagger\colon\Omega_M\rightarrow
\Omega_M$ the corresponding restriction maps. Evidently,
$(d_M^\dagger)^2=0$ and
${\star}^2=(-)^{(n-\partial)\partial}$.

The analog of the Laplace operator $\lM\colon\WM\rightarrow
\WM$ can be defined by the classical formula
\begin{equation}
\lM=\dM d^\dagger_M+d^\dagger_M \dM=(\dM+ d^\dagger_M)^2
\label{Lap}
\end{equation}
Evidently
\begin{align}
\lM \dM&=\dM\lM\\
\lM d^\dagger_M&=d^\dagger_M\lM\\
\lM \hM&=\hM\lM.
\end{align}
The Laplacian $\lM$ is grade-preserving. If $f\in\cal V$ then
\begin{equation}
\lM(f)=d^\dagger_M \dM(f)=-\sum_{i=1}^n X_i^2(f)
\end{equation}

Now we shall prove that, under certain additional assumptions,
it is possible to define the integration and scalar product in
$\WM$, such that  $d^\dagger_M$  and  $\dM$  become  mutually
formally adjoint.

As first let us  suppose  that $\cal  B$  can  be
realized as a dense *-subalgebra of a $C^*$-algebra $\widehat{\cal
B}$.

Secondly,
let us assume that there exists a faithful state
$\nu$ on
$\widehat{\cal B}$ such that
\begin{align}
\nu X_i(b)&=0\label{nuX}\\
\sum_k \nu(b_k)\otimes c_k&=\nu(b)\otimes 1\label{nu-cov}
\end{align}
for each $b\in\cal B$ and $i\in\{1,\ldots,n\}$, where
$\Sum_k b_k\otimes c_k=F(b)$.

Let $\Int_P\colon\Omega_P\rightarrow \BC$ be a
linear functional given by
\begin{equation*}
\int_P\Bigl(b\otimes (\e_1\wedge\ldots \wedge\e_n)\otimes W\Bigr)
=\nu(b)\qquad
\int_P\Bigl(\Omega_P^k\Bigr)=\{0\},
\end{equation*}
for $k<n+m$. Here $W\in L^{\wedge m}$ is a volume element
($m=\dim L=n(n-1)/2$).
\begin{lem} We have
\begin{align}
\sum_k\biggl(\int_P w_k\biggr)\otimes c_k&=
\biggl(\int_Pw\biggr)\otimes 1\label{int-1}\\
\int_Pd^\wedge(w)&=0\label{int-2}
\end{align}
for each $w\in\Omega_P$, where $\Sum_k w_k\otimes c_k=F^\wedge(w)$.
\end{lem}
\begin{pf}
Formula \eqref{int-1}   directly   follows   from
\eqref{nu-cov} and from the definition of $\Int_P$.
Formula \eqref{int-2} is non-trivial only if
$\partial w=n+m-1$. In this case we can write
$$ w=\sum_{i=1}^n b_i\otimes (\e_1\wedge\ldots
\wedge\e_{i-1}\wedge\e_{i+1}\wedge\ldots\wedge
\e_n)\otimes W +\sum_jq_j\otimes (\e_1\wedge\ldots
\wedge\e_n) \otimes\eta_j $$
where $\eta_j\in (L^{\wedge})^{m-1}$ (and $b_i,q_j\in\cal{B}$).
Applying \eqref{nuX}--\eqref{nu-cov}, the definition
of $d^\wedge$ and the fact that
$d(L^{\wedge m-1})=\{0\}$ we obtain
\begin{multline*}
\int_P d^\wedge(w)=-\sum_{i=1}^n (-)^i
\int_P\Bigl(X_i(b_i)\otimes
(\e_1\wedge\ldots \wedge\e_n) \otimes W\Bigr)\\
=-\sum_{i=1}^n (-)^i \nu\bigl(X_i(b_i)\bigr)=0. \qed
\end{multline*}
\renewcommand{\qed}{}
\end{pf}

Let $\intM\colon\hor_P\rightarrow\BC$ be a
linear functional given by
$$\intM\varphi=\int_P(\varphi\otimes W).$$
\begin{lem}
We have
\begin{equation}\label{int-3}
\intM\nabla(\varphi)=0
\end{equation}
for each $\varphi\in\hor_P$.
\end{lem}
\begin{pf}
Applying \eqref{int-2} and \eqref{d-hor}--\eqref{d-L}, and
the above definition we find
$$ \intM\nabla(\varphi)=\int_P\nabla(\varphi)\otimes W=
\int_P d^\wedge(\varphi)W=\int_P d^\wedge(\varphi W)=0.\qed$$
\renewcommand{\qed}{}
\end{pf}

The formula
\begin{equation}
(\psi, \varphi)=\intM \bar{\psi}{\star}(\varphi)
\end{equation}
where the bar denotes the standard *-operation,
determines a (strictly positive) scalar product in
$\hor_P$  (and
in particular in $\WM$). We have
\begin{equation}
(\nabla\psi,\varphi)=(\psi,\nabla^\dagger\varphi)
\end{equation}
for  each  $\varphi,\psi\in\hor_P$. In particular, $\dM$ and
$d_M^\dagger$ are formally adjoint in $\WM$ and the operator
$\lM$ is symmetric and positive. The map ${\star}$ is isometric.

\section{Associated Spinor Bundles}

Let $\S$ be the space of algebraic spinors for $\cl_n$.
Explicitly \cite{E}, $\S$ can be constructed as  follows. For simplicity
it will be assumed that $n$ is even. Then
$\BC_n=V_+\oplus V_-$ where  $V_\pm$ are (naturally mutually dual)
isotropic
subspaces spanned  by vectors  $\e_{2j-1}\pm i\e_{2j}$ while
$j\in\{1,\ldots,r\}$ and $2r=n$. By definition
$$\S=V_+^\wedge. $$
The space $\S$ is an irreducible left $\cl_n$-module.  Elements  of
the generating space $\BC_n$ of $\cl_n$ act on $\S$ as follows
$$ x\xi=x_+\wedge\xi+x_-\JJ\xi$$
where  $x_\pm\in  V_\pm$  are such that
$$x=\frac{x_-+x_+}{\sqrt{2}}$$
and $\JJ$ is the contraction map.

Let $\bigl\{\,\s_I\mid I\subseteq\{1,\ldots,r\}\,\bigr\}$
be a basis in  $\S$  given
by
$$ \s_I=\prod_{j\in I}\left(\frac{\e_{2j-1}+i\e_{2j}}{\sqrt{2}}
\right).$$
It is natural to introduce a scalar
product $(,)$ in $\S$ by requiring that $\s_I$ are
orthonormal  spinors.  Then  the  action  of  $\cl_n$  on  $\S$  is
hermitian, in the  sense  that $(Q\xi,\zeta)=(\xi,Q^*\zeta)$  for
each     $\xi,\zeta\in     \S$     and      $Q\in\cl_n$,      where
$*\colon\cl_n\rightarrow\cl_n$  is  the
antimultiplicative antilinear  involution  extending  the  complex
conjugation in $\BC_n$.
We shall denote by $\gamma_i\colon\S\rightarrow\S$ operators representing
the basis vectors $e_i$.

The spinorial conjugation operator $C\colon\S
\rightarrow\S$ is given by
$$ C(\psi)=(-)^{k(k+1)/2}{\star}_{\S}(\bar{\psi}),$$
where $\psi\in\S^k$ and ${\star}_\S$ is the
corresponding Hodge *-operator.
The bar denotes the complex conjugation in $\S$, defined by
requiring reality of basis vectors $\s_I$. We have
$$C\gamma_i+\gamma_iC=0$$
for each $i\in\{1,\ldots,n\}$. Further, $C$ is isometric and
$$C^2=(-)^{r(r+1)/2}.$$

Restricted on $\Spin (n)$ the action of $\cl_n$ becomes a unitary
representation $U$. We have
\begin{equation}\label{spin-dec}
\S=\S_+\oplus \S_-
\end{equation}
where $\S_+$ and $\S_-$ are subspaces corresponding to even and  odd
multivectors. These subspaces  are  invariant  under $U$,  moreover
$U\rest \S_\pm$ are irreducible. The operator $C$ intertwines $U$.

We shall also interpret the representation $U$ as a
right comodule structure
map $U\colon \S\rightarrow \S\otimes{\cal S}_n$.

In the case $n=2r+1$ similar algebraic constructions can be performed.

Now, the associated spinor bundle will be defined.
Let $\cal E$ be the space of all linear maps $\varsigma\colon
\S^*\rightarrow{\cal B}$ satisfying
\begin{equation}
F\varsigma=(\varsigma\otimes \id)U_c
\end{equation}
(intertwiners between the contragradient representation
$U_c\colon \S^*\rightarrow \S^*\otimes {\cal S}_n$ and $F$).
The space $\cal E$ is a  $\cal V$-bimodule, in  a  natural
manner.
The elements of $\cal E$ are interpretable as ``smooth
sections'' of the
``associated spinor bundle'' $E$.

Similarly, intertwiners $\varsigma\colon
\S^*\rightarrow\hor_P$ between $U_c$ and $F^\sstar$ constitute a
graded $\WM$-bimodule $\EW$  the  elements  of
which are interpretable as $E$-valued differential forms on
$M$. We have
$$\EW=\sideset{}{^\oplus}
\sum_{k=1}^n {\cal E}_\Omega^k$$
with ${\cal E}_\Omega^0={\cal E}$.

Equivalently, $\EW$  can  be
described as
the space of elements $\psi\in \hor_P\otimes  \S$  invariant  under
the product of actions $F^\sstar$  and  $U$. The connection between
the second and the first definition is given by
$\psi\leftrightarrow\varsigma$, with
\begin{gather*}
\varsigma(\zeta)=
(\id\otimes \zeta)(\psi)\\
\psi=\sum_I\varsigma(\s_I^*)\otimes \s_I
\end{gather*}
where $\s_I^*\in \S^*$ are vectors of the
corresponding biorthogonal basis.

The   decomposition \eqref{spin-dec}   induces a bimodule
splitting
\begin{equation}
\EW={\cal E}_\Omega^+\oplus{\cal E}_\Omega^-.
\end{equation}
In particular,
$ {\cal E}={\cal E}^+\oplus{\cal E}^-. $
Every $\psi\in\EW$ can be written in the form
\begin{equation}
\psi=\sum_I \psi_I\otimes \s_I
\end{equation}
where $\psi_I\in\hor_P$ and
\begin{equation}
F^\sstar(\psi_I)=\sum_J\psi_J\otimes U_{JI}^*.
\end{equation}

The space $\cal E$ possesses a natural *-structure, given by a
``charge conjugation'' $c_E\colon{\cal E}\rightarrow\cal{E}$, where
$$c_E(\psi)=\sum_I\bar{\psi}_I\otimes C(\s_I)$$
(the same formula determines a *-stucture on the
$\Omega_M$-bimodule $\cal{E}_\Omega$).

The corresponding Hodge operator
$\hE\colon{\cal  E}_\Omega\rightarrow{\cal  E}_\Omega$ is given by
\begin{equation}
\hE(\psi)=\sum_I {\star}(\psi_I)\otimes \s_I.
\end{equation}

By construction,
\begin{align*}
c_E^2&=(-)^{r(r+1)/2}\\
{\star}_E^2&=(-)^{\partial(n-\partial)}.
\end{align*}

The formula
\begin{equation}
\langle\psi,\varphi\rangle=\sum_I(\psi_I,\varphi_I)
\end{equation}
defines  a  (strictly  positive)  scalar  product   on
$\EW$.
The spaces ${\cal E}_\Omega^k$ are mutually  orthogonal. The operator
${\star}_E$ is isometric. It is worth noticing that $c_E$ is generally not
isometric (if the state $\nu$ possesses non-trivial modular properties).

Let $\nE,\nabla^\dagger_{\! E}\colon\EW
\rightarrow\EW$ be operators given by
\begin{align}
\nE(\psi)&=(\nabla^{\phantom{\dagger}}\!\otimes\id)\rest\cal{E}_\Omega\\
\nabla^\dagger_{\! E}(\psi)&=(\nabla^\dagger\!\otimes
\id)\rest\cal{E}_\Omega
\end{align}
In particular,
\begin{equation}
\nabla^\dagger_{\! E}(\psi)=-(-)^{nk-n}\hE\nE\hE(\psi)
\end{equation}
for each $\psi\in{\cal E}_\Omega^k$.

Operators $\nE$ and $\nabla^\dagger_{\! E}$ are mutually  formally
adjoint. The map  $\nE$  plays  the  role  of  the  spinorial
covariant derivative.

The corresponding Laplacian can be introduced by the classical
expression
\begin{equation}
\lE=\nE \nabla^\dagger_{\! E}+\nabla^\dagger_{\! E}\nE.
\end{equation}
We have
\begin{align*}
\langle\lE\psi,\varphi\rangle&=\langle\psi,\lE\varphi\rangle\\
\lE\hE&=\hE\lE.
\end{align*}
The operator $\lE$ is grade-preserving and positive.
In the ``pure spinor'' space $\cal E$ we have $\lE=
\nabla^\dagger_{\! E}\nE$. Explicitly,
\begin{equation}
\lE(\psi)=-\sum_{iI}X_i^2(\psi_I)\otimes \s_I
\end{equation}
for each $\psi\in\cal E$.

The map $c_E$ commutes with $\nabla_{\! E}$,
$\nabla_{\! E}^\dagger$, $\Delta_E$ and ${\star}_E$.

Now we shall introduce an analog of the Dirac operator. Let
$D_E\colon{\cal E}\rightarrow{\cal E}$ be a linear map defined by
\begin{equation}
D_E(\psi)=-i\sum_{iI}X_i(\psi_I)\otimes \gamma_i\s_I.
\end{equation}
By construction, $D_E$  ``mixes''  ${\cal  E}^+$  and  ${\cal  E}^-$.
Further,
\begin{align}
\langle D_E\psi,\varphi\rangle&=\langle\psi,D_E\varphi\rangle\\
D_Ec_E&=c_ED_E.
\end{align}

Let us consider ``the curvature operator'' $\varrho_E\colon{\cal E}
\rightarrow{\cal E}$ given by
\begin{equation}\label{curv-def}
\varrho_E(\psi)=-\frac{1}{2}\sum_{ijIJ}\psi_I\varrho_{ij}\bigl(
U_{JI}\bigr)\otimes\gamma_i\gamma_j\s_J
\end{equation}
where $\varrho_{ji}=-\varrho_{ij}\colon{\cal S}_n\rightarrow
{\cal B}$ are components of the curvature,
$$\varrho(a)=\frac{1}{2}\sum_{ij}\varrho_{ij}(a)\otimes
(\e_i\wedge\e_j).$$

\begin{lem} We have
\begin{equation}\label{D2=LR}
D_E^2=\lE+\varrho_E.
\end{equation}
\end{lem}
\begin{pf}
A direct calculation gives
\begin{equation*}
\begin{split}
D_E^2(\psi)&=-\sum_{ijI}X_iX_j(\psi_I)\otimes\gamma_i\gamma_j(
\s_I)\\
&=-\frac{1}{2}\sum_{ijI}[X_i,X_j](\psi_I)\otimes\gamma_i\gamma_j(
\s_I)-\sum_{iI}X_i^2(\psi_I)\otimes \s_I\\
&=\frac{1}{2}\sum_{ijIJ}\psi_{J}\varrho_{ij}(U_{JI}^*)
\otimes\gamma_i\gamma_j\s_I+
\lE(\psi)
=\varrho_E(\psi)+\lE(\psi).\qed
\end{split}
\end{equation*}
\renewcommand{\qed}{}
\end{pf}

Actually, the operator $\varrho_E$ can be further simplified so that
a complete analogy with the classical formalism \cite{Li}
holds.

The components of the curvature can be written as
\begin{equation}\label{comp}
\varrho_{ij}(a)=\frac{1}{2}\sum_{kl}\varrho_{ijkl}
e_{kl}\pi(a),
\end{equation}
where $\varrho_{ijkl}$
belong to the center of $\cal{B}$ and
$e_{kl}=-e_{lk}$ are canonical generators of $\so(n)=L^*$.

\begin{lem}We have
\begin{equation}
\varrho_E(\psi)=\frac{1}{4}\r\psi
\end{equation}
for each $\psi\in\cal{E}$, where
\begin{equation}
\r=\sum_{ij}\varrho_{ijij}.
\end{equation}
\end{lem}

The element $\r$ is interpretable as {\it the scalar curvature}
of the Levi-Civita connection.
\begin{pf}
The operator $\varrho_E$ can be transformed in the following way
\begin{equation*}
\begin{split}
\varrho_E(\psi)&=-\frac{1}{2}\sum_{ijIJ}\psi_I\varrho_{ij}\bigl(
U_{JI}\bigr)\otimes\gamma_i\gamma_j\s_J
=-\frac{1}{8}\sum_I\sum_{ijkl}\psi_I\varrho_{ijkl}\otimes
\gamma_i\gamma_j\gamma_k\gamma_l\s_I\\
&=\frac{1}{4}\sum_I\sum_{ij}\psi_I\varrho_{ijij}\otimes\s_I=
\frac{1}{4}\r\psi.
\end{split}
\end{equation*}
Here, we have applied the correspondence
$$e_{kl}\,\leftrightarrow\,\frac{1}{4}[\gamma_k,\gamma_l]$$
(the action of $\so(n)$ in the spinor space), as well as the
following {\it symmetry properties}
\begin{gather}
-\varrho_{ijkl}=\varrho_{jikl}=
\varrho_{ijlk}\\
\varrho_{ijkl}+\varrho_{iklj}+\varrho_{iljk}=0\\
\varrho_{ijkl}=\varrho_{klij}.
\end{gather}
They can be derived
essentially in the same manner as in classical geometry.
\end{pf}

The analog of the Clifford bundle can be  introduced  as  follows.
Let us consider a *-homomorphism
$\widehat{u}\colon\cl_n\rightarrow\cl_n\otimes
{\cal S}_{n}$ describing the natural (adjoint)
action of $\Spin(n)$ on $\cl_n$.
Let $\cl_M\subseteq{\cal B}\otimes
\cl_n$ be the fixed-point *-subalgebra for the  product  of  actions
$F$ and $\widehat{u}$. The elements of $\cl_M$ are  interpretable  as
``smooth sections'' of the corresponding  ``associated Clifford
bundle'' for $M$. By construction, $\cal V$ can be viewed as a
subalgebra of $\cl_M$. The space $\cal E$
is a left  $\cl_M$-module,  in  a  natural  manner.
\section{Example}
A large class of ``non-commutative spin manifolds''
can be constructed starting from appropriate
actions of $\Spin (n+1)$ on quantum spaces and  then  restricting
this action on $\Spin (n)$ (viewed as a subgroup of $\Spin (n+1)$,
in a natural manner).
At  the  dual  level,  this
inclusion of groups is described by an
epimorphism $j_n\colon {\cal S}_{n+1}\rightarrow  {\cal S}_n$  of
*-Hopf algebras (the restriction map).

Let ${\cal B}$ be a *-algebra. Let us assume that
$F\colon{\cal B}\rightarrow
{\cal B}\otimes{\cal S}_{n+1}$ is a *-homomorphism  describing  an
action of $\Spin (n+1)$ by *-automorphisms of ${\cal B}$
and satisfying the ``freeness'' condition. In  other
words for each $a\in {\cal S}_{n+1}$ there exist
$q_i,b_i\in{\cal B}$ such that \eqref{free} holds.

The formula
$$\iota(x)=\left(\begin{array}{cc}
0&x\\
-x^\top&0\end{array}\right)$$
defines an
embedding $\iota\colon\BC_n\rightarrow \so(n+1)$.

Let $F^\sharp\colon{\cal B}\rightarrow{\cal B}\otimes
{\cal S}_n$ be  the
restriction  of  $F$   on   $\Spin (n)$.   In   other   words
$$F^\sharp=(\id\otimes j_n)F.$$
Let ${\cal V}\subseteq{\cal B}$ be the fixed-point  subalgebra  of
${\cal B}$ for the action $F^\sharp$, and let
$i\colon{\cal V}\hookrightarrow {\cal B}$ be the inclusion map. Then
$P=({\cal B},i,F^\sharp)$ is a quantum
principal $\Spin (n)$-bundle (over the quantum space $M$ described
by ${\cal V}$).

The bundle $P$ possesses a canonical (integrable)  frame  structure.
For every $i\in\{1,\ldots,n\}$,     let     $X_i\colon{\cal
B}\rightarrow {\cal B}$ be a derivation corresponding
to $\iota(\e_i)$, via the infinitezimalization of  $F^\sharp$.  Let
$\partial_i=(X_i\rest{\cal V})\colon{\cal V}\rightarrow{\cal B}$ be
the corresponding restrictions.

Then $\tau=(\partial_1,\ldots,
\partial_n)$ is an integrable frame structure with
$\widehat{\tau}=(X_1,\ldots,X_n)$.  The  space   $M$   is   inherently
``curved''. The components $\varrho_{ij}$ of  the  curvature  take
values in scalars. Explicitly,
\begin{equation}
\varrho_{ij}(a)=-\bigl\{[\iota(\e_i),\iota(\e_j)](a)\bigr\}_\epsilon 1
\end{equation}
for each $a\in{\cal S}_n$.

As a concrete illustration, let us consider bundles based on
Cuntz algebras \cite{cun}. By definition the Cuntz algebra  ${\cal
C}_d$ (where $d\geq2$ is an integer) is a *-algebra  generated
by elements $\psi_1,\ldots,\psi_d$ and relations
\begin{gather}
\psi_i^*\psi_j=\delta_{ij}1\label{rel1}\\
\sum_{i=1}^d\psi_i\psi_i^*=1.
\end{gather}

Let us assume that $v\colon \BC_d\rightarrow\BC_d
\otimes{\cal S}_{n+1}$  is  a  faithful  and  irreducible  unitary
representation of $\Spin (n+1)$.
This map is uniquely extendible to
a *-homomorphism
$F\colon{\cal C}_d\rightarrow{\cal C}_d\otimes {\cal S}_{n+1}$
(describing the action of $\Spin (n+1)$ by automorphisms of ${\cal C}_d$).

Let us prove  that  the  action  $F$ satisfies the ``freeness''
condition. Matrix   elements   $v_{ij}$
generate the whole algebra ${\cal S}_{n+1}$. For this reason it is
sufficient to check that  equations  of  the  form \eqref{free}
hold  for elements $a=v_{ij}$. We have
$$F(\psi_i)=\sum_{k=1}^d\psi_k\otimes v_{ki}.$$
Multiplying both sides of the above equality by $\psi_j^*$ on  the
left and using \eqref{rel1} we obtain
$$ \psi_j^*F(\psi_i)=1\otimes v_{ji}. $$
Consequently, the freeness condition holds.

\section{Concluding Remarks}
Logically, constructions of the Laplacian $\lM$
and the  Hodge  operator $\hM$ (and considerations with the curvature)
do not require
$\Spin (n)$-bundles. Everything
can be performed starting from principal $\SO(n)$-bundles
(and integrable frame structures on them).

The bundle of ``orthonormal oriented frames'' can be extracted from
the spin-bundle $P$  as  follows.
Let ${\cal O}\subseteq{\cal B}$ be  a
*-subalgebra    consisting    of    elements    $b$     satisfying
$F_g(b)=F_{-g}(b)$. By construction, $i(\cal V)\subseteq{\cal O}$
and $F({\cal  O})\subseteq{\cal  O}\otimes{\cal  A}_n$.  In  other
words, the   action   of
$\Spin (n)$ can be ``projected'' to the action of $\SO(n)$ on  ${\cal
O}$.  In  such  a  way  we  obtain  a   quantum   principal
$\SO(n)$-bundle $P^\sstar$ over  $M$.
Moreover, $X_i({\cal O})
\subseteq{\cal O}$  and  the  frame  structure  on  $P$  induces
a  frame
structure on $P^\sstar$. The bundle $P^\sstar$ is
interpretable  as  consisting  of ``oriented orthonormal frames''
of $M$. At the level of spaces the inclusion ${\cal O}
\hookrightarrow{\cal B}$ becomes the ``covering projection''
of  $P$  onto
$P^\sstar$, corresponding to  the  standard  interpretation  of  spinor
structures.

The  number  $n$  (playing  the  role  of  the
dimension of $M$) is not  fixed.  The  same  quantum  space  could
possess frame structures with different ``dimensions''.  Moreover,
even if $M$ is a classical compact  manifold  then  the
class of frame structures described above includes
nonstandard geometrical  objects
based on foliations of $M$.

It is interesting to think about possible physical applications of the
presented theory, starting from the idea that space-time is interpretable
as a non-commutative (completely pointless) spin manifold. In this conceptual
framework it is natural to consider field theories including ``spinor''
and ``metric'' fields, as well as the corresponding ``gauge fields''.
For example, classical Einstein-Dirac type theories
\cite{CruE1,CruE2}
include all these concepts in a natural way. Counterparts of such theories
(on non-commutative spin manifolds) are in principle constructable
following the ideas of classical geometrical formulations.

\end{document}